\documentclass[12pt,a4paper]{article}
\usepackage{amsmath,amsthm,amssymb,latexsym,epic,bezier}

\newtheorem{theorem}{Theorem}
\newtheorem{lemma}[theorem]{Lemma}
\newtheorem{corollary}[theorem]{Corollary}
\newtheorem{proposition}[theorem]{Proposition}

\newtheorem{remark}[theorem]{Remark}
\newtheorem{question}[theorem]{Question}
\usepackage[all]{xy}
\usepackage{enumerate}
\numberwithin{equation}{section}

\renewcommand{\phi}{\varphi}

\begin{document}

\title{On Kiselman's semigroup}
\author{Ganna Kudryavtseva and Volodymyr Mazorchuk}
\date{}
\maketitle

\begin{abstract}
We study the algebraic properties of the series
$\mathrm{K}_n$  of semigroups, which is inspired by \cite{Ki}
and has origins in convexity theory.
In particular, we describe Green's relations on $\mathrm{K}_n$, 
prove that there exists a faithful representation 
of $\mathrm{K}_n$ by $n\times n$ matrices with non-negative 
integer coefficients (and even explicitly construct such a 
representation), and prove that $\mathrm{K}_n$ does not admit a 
faithful representation by matrices of smaller size. 
We also  describe the maximal nilpotent subsemigroups in
$\mathrm{K}_n$, all isolated and completely isolated subsemigroups,
all automorphisms and anti-automorphisms of $\mathrm{K}_n$.
Finally, we explicitly construct all irreducible representations 
of $\mathrm{K}_n$ over any field and describe primitive 
idempotents in the semigroup algebra (which we prove is basic).
\end{abstract}

\section{Introduction}\label{s1}

Let $E$ be a real vector space and $\mathrm{Func}(E)$ be the set of all
functions on $E$ with values in the extended real line
$\mathbb{R}\cup\{-\infty,+\infty\}$. In convexity theory there
appear three natural operators on $\mathrm{Func}(E)$, namely the
operator $c$ of taking the convex hull of a function, the operator
$l$ of taking the largest lower semicontinuous minorant of the function,
and the operator $m$ defined via
\begin{displaymath}
m(f)(x)=
\begin{cases}
f(x),& \text{if $f$ is everywhere }>-\infty;\\
-\infty,& \text{otherwise}.
\end{cases}
\end{displaymath}
The operators $c,l,m$ generate a monoid, $G(E)$, with repsect to the
usual composition. In \cite{Ki} it was shown that this monoid 
consists of $18$ elements and has the following presentation
(as a monoid):
\begin{multline}\label{eq1}
G(E)=\langle c,l,m:c^2=c,\,\, l^2=l,\,\, m^2=m,\\
clc=lcl=lc,\,\, cmc=mcm=mc,\,\, lml=mlm=ml\rangle.
\end{multline}
Furthermore, the paper \cite{Ki} also contains a detailed study of
the algebraic structure of $G(E)$ and gives a faithful representation
of $G(E)$ by $3\times 3$ matrices with non-negative integer coefficients.

There is a fairly straightforward way to generalize \eqref{eq1}. 
Let $n$ be a positive integer. Denote by $\mathrm{K}_n$ the monoid
defined via the following presentation:
\begin{multline}\label{eq2}
\mathrm{K}_n=\langle a_1,\dots,a_n:\,\, a_i^2=a_i,\,i=1,\dots,n;\\
a_ia_ja_i=a_ja_ia_j=a_ja_i,\, 1\leq i<j\leq n\rangle.
\end{multline}
We will call $\mathrm{K}_n$ {\em Kiselman's semigroup} after the author 
of \cite{Ki}. Obviously, we have $G(E)\cong \mathrm{K}_3$. 
The generalization \eqref{eq2} was proposed 
by O.~Ganyushkin and the second author in 2002 (unpublished). 
In \cite{Go} several results on the structure of
$\mathrm{K}_n$ were announced. Unfortunately, the proofs have never
appeared. So, we have decided to study $\mathrm{K}_n$ independently.
In the present paper we prove all the results announced in
\cite{Go}, in particular, we describe Green's relations on 
$\mathrm{K}_n$ (Section~\ref{s7}), 
prove that there exists a faithful representation 
of $\mathrm{K}_n$ by $n\times n$ matrices with non-negative 
integer coefficients (and even explicitly construct such a 
representation), and  prove that $\mathrm{K}_n$ does not admit a 
faithful representation by matrices of smaller size
(Subsection~\ref{s11.1}).  We also obtain some additional results, in
particular, we describe the maximal nilpotent subsemigroups in
$\mathrm{K}_n$ (Section~\ref{s8}), all isolated and 
completely isolated subsemigroups (Section~\ref{s9}),
all automorphisms of $\mathrm{K}_n$ and all anti-automorphisms 
of $\mathrm{K}_n$ (Section~\ref{s6}). We also explicitly construct
all irreducible representations of $\mathrm{K}_n$ over any field
and describe the primitive idempotents in the semigroup algebra
(Subsection~\ref{s11.2}). We are convinced that $\mathrm{K}_n$ is 
a very beautiful  combinatorial objects and might have a lot of 
further interesting combinatorial properties and applications.

\vspace{5mm}

\noindent
{\bf Acknowledgments.} The paper was written during the visit of the 
first author to Uppsala University, which was supported by the 
Swedish Institute. The financial support of the Swedish Institute and 
the hospitality of Uppsala University are gratefully acknowledged. For 
the second author the research was partially supported by the Swedish 
Research Council. 
\vspace{0.5cm}

\section{Finiteness of $\mathrm{K}_n$}\label{s2}

We will denote by $e$ the unit element in $\mathrm{K}_n$.
For a finite alphabet, $\mathcal{A}$, we denote by
$\mathrm{W}(\mathcal{A})$ the set of all  {\em finite} words over 
this alphabet, including the empty word (with respect to the usual
operation of concatenation of words this is the same as
the free monoid, generated by $\mathcal{A}$, which is sometimes
denoted by $\mathcal{A}^*$). Let $\mathfrak{l}:\mathrm{W}(\mathcal{A})\to\mathbb{N}\cup\{0\}$
denote the {\em length function}.

\begin{lemma}\label{lemma1}
\begin{enumerate}[(i)]
\item\label{lemma1.1} Let $i\in\{1,\dots,n\}$ and 
$w\in \mathrm{W}(\{a_1,\dots,a_{i-1}\})$. Then we have
$a_iwa_i=a_iw$ in $\mathrm{K}_n$.
\item\label{lemma1.2} Let $i\in\{1,\dots,n\}$ and 
$w\in \mathrm{W}(\{a_{i+1},\dots,a_{n}\})$. Then we have
$a_iwa_i=wa_i$ in $\mathrm{K}_n$.
\end{enumerate}
\end{lemma}

\begin{proof}
We prove \eqref{lemma1.1}. The statement \eqref{lemma1.2} is
proved by similar arguments. We proceed by induction on
$\mathfrak{l}(w)$. If $\mathfrak{l}(w)=0$ or
$\mathfrak{l}(w)=1$, the statement follows directly from
the presentation \eqref{eq2}. Assume now that 
$\mathfrak{l}(w)>1$ and write $w=w'a_j$ for some $j<i$.
Then $w'\in \mathrm{W}(\{a_1,\dots,a_{i-1}\})$ and
$\mathfrak{l}(w')=\mathfrak{l}(w)-1$.
We have
\begin{multline*}
a_iwa_i=a_iw'a_ja_i=(a_iw')a_ja_i=
\text{(by the inductive assumption)}=\\
=(a_iw'a_i)a_ja_i=a_iw'a_ia_ja_i=\text{(by \eqref{eq2})}=
a_iw'a_ia_j=(a_iw'a_i)a_j=\\=
\text{(by the inductive assumption)}=(a_iw')a_j=a_iw'a_j=a_iw.
\end{multline*}
\end{proof}

Define the function $L:\mathbb{N}\to \mathbb{N}$ as follows:
\begin{displaymath}
L(n)=
\begin{cases}
2^{k+1}-2,& n=2k;\\
3\cdot 2^{k}-2,& n=2k+1.
\end{cases}
\end{displaymath}

\begin{corollary}\label{cnew1}
Let $\alpha\in \mathrm{K}_n$, $\alpha\neq e$, and let 
$w\in \mathrm{W}(\{a_1,\dots,a_{n}\})$ be a word of the 
shortest possible length such that $\alpha=w$ in $\mathrm{K}_n$.
Then we have the following:
\begin{enumerate}[(i)]
\item\label{cnew1.1} For $i\leq \lceil\frac{n}{2}\rceil$
the letter $a_i$ occurs in $w$ at most $2^{i-1}$ times.
\item\label{cnew1.2} For $i\geq \lceil\frac{n+1}{2}\rceil$
the letter $a_i$ occurs in $w$ at most $2^{n-i}$ times.
\item\label{cnew1.3} $\mathfrak{l}(w)\leq L(n)$.
\end{enumerate}
\end{corollary}

\begin{proof}
We prove \eqref{cnew1.1} by induction on $i$. If the letter 
$a_1$ occurs in $w$ more than once, the word $w$ can be 
reduced (shortened) using Lemma~\ref{lemma1}\eqref{lemma1.2}. 
This gives us the
basis of the induction. Let $1<i\leq \lceil\frac{n}{2}\rceil$. 
From the inductive 
assumption we obtain that the total number of occurrences of 
the letters $a_1,\dots,a_{i-1}$ in $w$ does not exceed
$2^{i-1}-1$. Hence we can write 
$w=w_1b_1w_2b_2w_3\dots w_{2^{i-1}-1}b_{2^{i-1}-1}w_{2^{i-1}}$,
where $b_j\in \{a_1,\dots, a_{i-1}\}$ and 
$w_j\in \mathrm{W}(\{a_{i},\dots,a_{n}\})$ for all appropriate
$j$. If $a_i$ occurs in some $w_j$ more than once, the word 
$w_j$ and hence $w$ can be reduced using 
Lemma~\ref{lemma1}\eqref{lemma1.2}. Hence the total number of
occurrences of $a_i$ in $w$ does not exceed $2^{i-1}$. This proves
\eqref{cnew1.1}. \eqref{cnew1.2} is proved by similar arguments.
\eqref{cnew1.3} follows from \eqref{cnew1.1} and \eqref{cnew1.2}
since for all $n=2k\in\mathbb{N}$ we have
\begin{displaymath}
L(n)=\sum_{i=1}^k 2^{i-1}+\sum_{i=k+1}^n 2^{n-i}
\end{displaymath}
and for all $n=2k+1\in\mathbb{N}$ we have
\begin{displaymath}
L(n)=\sum_{i=1}^{k+1} 2^{i-1}+\sum_{i=k+2}^n 2^{n-i}.
\end{displaymath}
\end{proof}

As an immediate corollary from the latter statement we have:

\begin{theorem}\label{theorem1}
The semigroup $\mathrm{K}_n$ is finite, moreover
\begin{displaymath}
|\mathrm{K}_n|\leq 1+n^{L(n)}.
\end{displaymath}
\end{theorem}

\begin{proof}
The semigroup $\mathrm{K}_n$ is generated by $n$ elements.
By Corollary~\ref{cnew1}\eqref{cnew1.3}, every element of 
$\mathrm{K}_n$, different from the unit element $e$, can be written 
as a product of at most $L(n)$ generators. Since all generators
are idempotents, repeating the last generator, occurring in such a
product, we conclude that every element of 
$\mathrm{K}_n$, different from the unit element $e$, can be written 
as a product of exactly $L(n)$ generators. The statement follows.
\end{proof}

\begin{question}\label{question1}
Can one give an explicit formula for $|\mathrm{K}_n|$?
\end{question}

\begin{remark}\label{remark1}
{\rm
In \cite{Go} a slightly more general family of semigroups is
considered: let $(I,<)$ be a partially ordered set. Define
\begin{displaymath}
\mathrm{K}_I=\langle a_i,\, i\in I:\,\, a_i^2=a_i,\,i\in I;\,\,
a_ia_ja_i=a_ja_ia_j=a_ja_i,\, i,j\in I,\, i<j\rangle.
\end{displaymath}
\cite[Theorem~2]{Go} states that $\mathrm{K}_I$ is finite if and
only if $I$ is finite and $<$ is linear. This is an immediate
consequence of Theorem~\ref{theorem1}. Indeed, Theorem~\ref{theorem1}
gives us the sufficiency. The necessity follows from the trivial
observation that for incomparable $i,j\in I$ the elements
$(a_ia_j)^k\in \mathrm{K}_I$, $k\in\mathbb{N}$, are obviously
different since there is no relation involving both $a_i$
and $a_j$.
}
\end{remark}

\section{The canonical form for elements of $\mathrm{K}_n$}\label{s3}

Let $\varphi:\mathrm{W}(\{a_1,\dots,a_{n}\})\to \mathrm{K}_n$
denote the canonical epimorphism. For $w\in 
\mathrm{W}(\{a_1,\dots,a_{n}\})$ set $\overline{w}=\{x\in
\mathrm{W}(\{a_1,\dots,a_{n}\}):\varphi(x)=\varphi(w)\}$.
If $w=a_{i_1}a_{i_2}\dots a_{i_k}\in \mathrm{W}(\{a_1,\dots,a_{n}\})$,
then by a {\em subword} of $w$ we will mean an element of
$\mathrm{W}(\{a_1,\dots,a_{n}\})$ of the form
$a_{i_s}a_{i_{s+1}}a_{i_{s+2}}\dots a_{i_t}$ for some 
$1\leq s\leq t\leq k$. By a {\em quasi-subword} of $w$ we will 
mean an element of $\mathrm{W}(\{a_1,\dots,a_{n}\})$ of the form
$a_{i_{l_1}}a_{i_{l_2}}a_{i_{l_3}}\dots a_{i_{l_t}}$ for some 
$1\leq l_1<l_2<l_3<\dots<l_t\leq k$ (including the empty 
quasi-subword). Each subword is, by definition, a quasi-subword.

The main result of this section is the following statement:

\begin{theorem}\label{theorem2}
Let $w\in \mathrm{W}(\{a_1,\dots,a_{n}\})$.
\begin{enumerate}[(i)]
\item\label{theorem2.1} 
The set $\overline{w}$ contains a unique element of the 
minimal possible length.
\item\label{theorem2.2}
$v\in \overline{w}$ has the minimal possible length if 
and only if the 
for each $i\in\{1,2,\dots,n\}$ the following condition is
satisfied: if $a_iua_i$ is a subword of $v$ then $u$ 
contains some $a_j$ with $j>i$ and some $a_k$ with $k<i$.
\end{enumerate}
\end{theorem}

The words $v\in \mathrm{W}(\{a_1,\dots,a_{n}\})$, satisfying the
condition of Theorem~\ref{theorem2}\eqref{theorem2.2}, will
be called {\em canonical}. If $w\in \mathrm{W}(\{a_1,\dots,a_{n}\})$
and $v\in\overline{w}$ is canonical, we will say that
$v$ is the {\em canonical form} of $w$.
By Theorem~\ref{theorem2}\eqref{theorem2.1}
the homomorphism $\varphi$ induces a bijection between the
set of all canonical words in $\mathrm{W}(\{a_1,\dots,a_{n}\})$
and the elements of $\mathrm{K}_n$. In particular, it
makes sense to speak about the {\em canonical form} of an
element from $\mathrm{K}_n$.

\begin{remark}\label{remark2}
{\rm
The statement of Theorem~\ref{theorem2}\eqref{theorem2.1} 
was announced in \cite[Theorem~1]{Go}.
}
\end{remark}

\begin{proof}
Define the binary relation $\to$ on $\mathrm{W}(\{a_1,\dots,a_{n}\})$
in the following way: for $w,v\in \mathrm{W}(\{a_1,\dots,a_{n}\})$
we set $w\to v$ if and only if there exists $i\in \{1,\dots,n\}$
such that $w=w_1a_iua_iw_2$ and either $v=w_1a_iuw_2$ and 
$u\in \mathrm{W}(\{a_1,\dots,a_{i-1}\})$, or  $v=w_1ua_iw_2$ and 
$u\in \mathrm{W}(\{a_{i+1},\dots,a_{n}\})$. From Lemma~\ref{lemma1}
we obtain that $w\to v$ implies $v\in\overline{w}$. Obviously,
$w\to v$ implies $\mathfrak{l}(v)=\mathfrak{l}(w)-1$, in particular,
any chain $v_1\to v_2\to v_3\to \dots$ in 
$\mathrm{W}(\{a_1,\dots,a_{n}\})$ terminates in a finite
number of steps. Denote by $\overset{*}{\to}$ the reflexive-transitive
closure of $\to$.

\begin{lemma}\label{lemma2}
For all $u,v,w\in \mathrm{W}(\{a_1,\dots,a_{n}\})$, such that
$u\neq v$, $w\to u$ and $w\to v$, there exists
$x\in \mathrm{W}(\{a_1,\dots,a_{n}\})$ such that
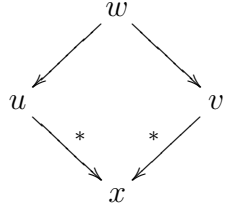
\begin{displaymath}
\xymatrix{
 &w\ar[dl]\ar[dr]& \\
u\ar[dr]^{*}&&v\ar[dl]_{*}\\
 &x&
}
\end{displaymath}
\end{lemma}

\begin{proof}
Both $u$ and $v$ are quasi-subwords of $w$ by the definition
of $\to$. $u$ is obtained from $w$ by deleting some $a_i$, and
$v$ is obtained from $w$ by deleting some $a_j$. If $i\neq j$,
from Lemma~\ref{lemma1} we obtain that we are allowed to 
delete the corresponding occurrence of $a_i$ in $v$ obtaining
some $x$ such that $v\to x$. Moreover, again applying 
Lemma~\ref{lemma1} we have that we are allowed to 
delete the corresponding occurrence of $a_j$ in $u$. Since
these operations obviously commute we will get the same result
$x$ and $u\to x$, as required. 

Now assume that $i=j$. By the definition of $\to$, the deletion
of $a_i$ involves two occurrences of $a_i$ in a word. If the
corresponding two pairs of $a_i$'s in $w$ do not intersect, then
the same argument as above works, implying that our deletion
operations commute. 

Without loss of generality, in the remaining cases we may assume
$w=a_i\alpha a_i\beta a_i$, where
$\alpha,\beta\in \mathrm{W}(\{a_1,\dots,a_{i-1},a_{i+1},\dots,a_n\})$.
If $u=v$, we can obviously take $x=u=v$. Hence we are left to 
deal with the following cases:
\begin{enumerate}
\item $u=\alpha a_i\beta a_i$, $v=a_i\alpha\beta a_i$. Because of
\eqref{eq2} this is possible if and only if $\alpha=e$, which gives
us $u=v$. This case was considered above.
\item $u=\alpha a_i\beta a_i$, $v=a_i\alpha a_i\beta$. 
In this case we can take $x=\alpha a_i\beta$ and obviously have
$u\to x$, $v\to x$.
\item $u=a_i\alpha a_i\beta $, $v=a_i\alpha\beta a_i$. Because of
\eqref{eq2} this is possible if and only if $\beta=e$, which gives
us $u=v$. This case was considered above.
\end{enumerate}
The statement of the lemma follows.
\end{proof}

The statement \eqref{theorem2.1} follows now from
Lemma~\ref{lemma2} and the Diamond Lemma (see e.g. \cite{Ne}).
The statement \eqref{theorem2.2} follows from the 
statement \eqref{theorem2.1} and the definition of the relation
$\to$. This completes the proof.
\end{proof}

From Corollary~\ref{cnew1}\eqref{cnew1.1} we know that for any
$w\in \mathrm{W}(\{a_1,\dots,a_{n}\})$ the length of the 
minimal representative in $\overline{w}$ does not
exceed $L(n)$. Now we can show that this bound is sharp.

\begin{corollary}\label{cnew2}
There exists $w\in \mathrm{W}(\{a_1,\dots,a_{n}\})$ such that
the length of the minimal representative in $\overline{w}$
equals $L(n)$.
\end{corollary}

\begin{proof}
Let $k=\lceil\frac{n}{2}\rceil$ and set $w_1=a_1a_n$,
$w_2=a_2a_{n-1}$,\dots, $w_{k-1}=a_{k-1}a_{n-k+2}$,
\begin{displaymath}
w_k=
\begin{cases}
a_{k}a_{n-k+1},& n\text{ is even};\\
a_k,& n\text{ is odd}.
\end{cases}
\end{displaymath}
Define the words $v_i$, $i=1,\dots,k$, recursively as follows:
$v_1=w_1$; if $v_i=w_{j_1}w_{j_2}\dots w_{j_s}$, then
$v_{i+1}=w_{i+1}w_{j_1}w_{i+1}w_{j_2}w_{i+1}\dots w_{i+1}w_{j_s}w_{i+1}$.
It follows immediately that $\mathfrak{l}(v_{k})=L(n)$ and
it is easy to see from the construction that $v_i$ is canonical
for every $i$. The claim follows.
\end{proof}

\section{Idempotents in $\mathrm{K}_n$}\label{s4}

Let $w\in \mathrm{W}(\{a_1,\dots,a_{n}\})$. Define the
{\em content} $\mathfrak{c}(w)$ of $w$ as the set of all
those $i\in\{1,\dots,n\}$ such that the letter $a_i$ appears
in $w$. In particular, $\mathfrak{c}(e)=\varnothing$ and
$\mathfrak{c}(a_i)=\{i\}$ for all $i=1,\dots,n$.
From \eqref{eq2} it follows immediately that
$\mathfrak{c}(v)=\mathfrak{c}(w)$ for every $v\in\overline{w}$,
in particular, one can speak of the {\em content} of an
element from $\mathrm{K}_n$. Furthermore, obviously
$\mathfrak{c}(wv)=\mathfrak{c}(w)\cup\mathfrak{c}(v)$ for
all $v,w\in \mathrm{W}(\{a_1,\dots,a_{n}\})$, which implies the
following statement:

\begin{lemma}\label{lemma4}
$\mathfrak{c}$ is an epimorphism from the semigroup $\mathrm{W}(\{a_1,\dots,a_{n}\})$ to the semigroup
$(2^{\{1,2,\dots,n\}},\cup)$. $\mathfrak{c}$ also induces 
an epimorphism from $\mathrm{K}_n$ to the semigroup
$(2^{\{1,2,\dots,n\}},\cup)$ (abusing notation we will
denote this epimorphism also by $\mathfrak{c}$).
\end{lemma}

Let $X\subset \{1,\dots,n\}$. If $X=\varnothing$, set
$e_{\varnothing}=e$. If $X\neq \varnothing$, let
$X=\{i_1,\dots,i_k\}$ such that $i_1>i_2>\dots>i_k$. Set
$e_X=a_{i_1}a_{i_2}\cdots a_{i_k}$.

\begin{proposition}\label{proposition2}
Each $e_X$ is an idempotent in $\mathrm{K}_n$ and every
idempotent in $\mathrm{K}_n$ has the form $e_X$ for some
$X\subset \{1,\dots,n\}$. In particular, the semigroup 
$\mathrm{K}_n$ contains $2^n$ idempotents.
\end{proposition}

\begin{proof}
As the word $a_{i_1}a_{i_2}\cdots a_{i_k}$ is canonical we have
$e_X\neq e_Y$ if $X\neq Y$. That 
$e_Xe_X=e_X$ follows immediately from Lemma~\ref{lemma1}\eqref{lemma1.1}.
Hence we have only to show that any idempotent in $\mathrm{K}_n$ 
has the form $e_X$ for some $X\subset \{1,\dots,n\}$. 
Let $x\in \mathrm{K}_n$ be an idempotent. Then
$x^k=x$ for all $k\in\mathbb{N}$ and the necessary statement
follows from the following lemma:

\begin{lemma}\label{lemma3}
Let $w\in \mathrm{W}(\{a_1,\dots,a_{n}\})$. Then
$w^{k}=e_{\mathfrak{c}(w)}$ for all $k\geq |\mathfrak{c}(w)|$.
\end{lemma}

\begin{proof}
Set $N=|\mathfrak{c}(w)|$.
Let $X\subset \{1,\dots,n\}$. From Lemma~\ref{lemma1}\eqref{lemma1.1} 
and the definition of $e_X$ it follows that $e_X a_i=e_X$ for every 
$i\in X$. Hence it is enough to show that
$w^{N}=e_{\mathfrak{c}(w)}$. For $i\in\{1,\dots,n\}$ denote by 
$\partial_i:\mathrm{W}(\{a_1,\dots,a_{n}\})\to
\mathrm{W}(\{a_1,\dots,a_{i-1},a_{i+1},\dots,a_n\})$ the operation of
deleting all occurrences of the letter $a_i$ in a word.
Let $\mathfrak{c}(w)=\{i_1,\dots,i_N\}$ and
$i_1>i_2>\dots>i_N$. Using Lemma~\ref{lemma1}\eqref{lemma1.1} 
we inductively compute:
\begin{multline}\label{eq3}
w^N=\underbrace{wwww\dots w}_{N \text{ times }}=
w\partial_{i_1}(w)\partial_{i_1}(w)\partial_{i_1}(w)\dots
\partial_{i_1}(w)=\\=
w\partial_{i_1}(w)\partial_{i_2}\partial_{i_1}(w)
\partial_{i_2}\partial_{i_1}(w)\dots
\partial_{i_2}\partial_{i_1}(w)=\dots=\\=
w\partial_{i_1}(w)\partial_{i_2}\partial_{i_1}(w)
\partial_{i_2}\partial_{i_3}\partial_{i_1}(w)\dots
(\partial_{i_{N-1}}\dots\partial_{i_2}\partial_{i_1})(w).
\end{multline}
For $j=1,\dots,N-1$ set
$w_j=\partial_{i_{j}}\dots\partial_{i_2}\partial_{i_1}(w)$.
Again, from the computation \eqref{eq3} and 
Lemma~\ref{lemma1}\eqref{lemma1.2}  we inductively derive:
\begin{multline}\label{eq4}
w^N=ww_1w_2\dots w_{N-1}=
\partial_{i_N}(w)\partial_{i_N}(w_1)
\partial_{i_N}(w_2)\dots\partial_{i_N}(w_{N-2})w_{N-1}=\\
=\partial_{i_{N-1}}\partial_{i_N}(w)\partial_{i_{N-1}}
\partial_{i_N}(w_1)\dots\partial_{i_{N-1}}\partial_{i_N}(w_{N-3})
\partial_{i_N}(w_{N-2})w_{N-1}=\dots=\\
=(\partial_{i_2}\dots\partial_{i_{N-1}}\partial_{i_N})(w)
(\partial_{i_3}\dots\partial_{i_{N-1}}\partial_{i_N})(w_1)
\dots\partial_{i_N}(w_{N-2})w_{N-1}.
\end{multline}
Now it is left to observe that 
\begin{multline*}
\mathfrak{c}((\partial_{i_2}\dots\partial_{i_{N-1}}
\partial_{i_N})(w))=\{i_1\},
\mathfrak{c}((\partial_{i_3}\dots\partial_{i_{N-1}}
\partial_{i_N})(w_1))= \{i_2\},\dots,\\
\mathfrak{c}(w_{N-1})= \{i_N\}.
\end{multline*}
Hence the product in the formula \eqref{eq4} results in the 
product $a_{i_1}a_{i_2}\dots a_{i_N}$,
which is equal to $e_{\mathfrak{c}(w)}$. Therefore 
$w^N=e_{\mathfrak{c}(w)}$ and the statement is proved.
\end{proof}
The statement of Proposition~\ref{proposition2} follows
immediately from Lemma~\ref{lemma3}.
\end{proof}

\begin{remark}\label{remark15}
{\rm 
It is easy to see that different idempotents in $\mathrm{K}_n$
do not commute. Furthermore, the set of all idempotents in
$\mathrm{K}_n$ is not a subsemigroup of $\mathrm{K}_n$, as it follows
from the next statement.
}
\end{remark}

\begin{proposition}\label{proposition70}
Let $X,Y\subset \{1,2,\dots,n\}$. Then the following conditions
are equivalent:
\begin{enumerate}[(a)]
\item\label{proposition70.1} $e_Xe_Y$ is an idempotent.
\item\label{proposition70.2} $e_Xe_Y=e_{X\cup Y}$.
\item\label{proposition70.3} For every $i\in X\setminus Y$
and every $j\in Y\setminus X$ we have $i>j$.
\end{enumerate}
\end{proposition}

\begin{proof}
The implication \eqref{proposition70.2}$\Rightarrow$\eqref{proposition70.1}
is obvious. By Lemma~\ref{lemma4} we have 
$\mathfrak{c}(e_Xe_Y)=X\cup Y$. At the same time $e_{X\cup Y}$ is the 
only idempotent of $\mathrm{K}_n$ with content
$X\cup Y$. The implication
\eqref{proposition70.1}$\Rightarrow$\eqref{proposition70.2} follows.

If $|X|=0$, the implication 
\eqref{proposition70.3}$\Rightarrow$\eqref{proposition70.2} is trivial.
Hence we may assume $|X|>0$.
We prove the implication 
\eqref{proposition70.3}$\Rightarrow$\eqref{proposition70.2} by induction
on $|Y|$. If $|Y|=0$, we have $e_Y=e$ and the claim is obvious.
Let $|Y|>0$ and $y$ be the minimal element of $Y$. Let $x$ be the
minimal element of $X$. If $x=y$, we have 
\begin{displaymath}
e_Xe_Y=e_{X\setminus\{x\}}a_ye_{Y\setminus\{y\}}a_y=
e_{X\setminus\{x\}}e_{Y\setminus\{y\}}a_y
\end{displaymath}
by Lemma~\ref{lemma1}\eqref{lemma1.2}. The sets 
$X\setminus\{x\}$ and $Y\setminus\{y\}$ still satisfy
\eqref{proposition70.3} and hence by induction we get
\begin{displaymath}
e_{X\setminus\{y\}}e_{Y\setminus\{y\}}a_y=
e_{(X\cup Y)\setminus\{y\}}a_y=e_{X\cup Y}.
\end{displaymath}

If $x\neq y$, then $x>y$ by \eqref{proposition70.3}. Hence
the sets $X$ and $Y\setminus\{y\}$ satisfy
\eqref{proposition70.3} and hence by induction we get
\begin{displaymath}
e_Xe_Y=e_{X}e_{Y\setminus\{y\}}a_y=
e_{(X\cup Y)\setminus\{y\}}a_y=e_{X\cup Y}.
\end{displaymath}
This proves the implication
\eqref{proposition70.3}$\Rightarrow$\eqref{proposition70.2}.

Finally, assume that \eqref{proposition70.3} is not satisfied.
Let $i\in X\setminus Y$ and $j\in Y\setminus X$ be such that $i<j$.
Then the letter $a_i$ occurs in $e_Xe_Y$ to the left of the 
letter $a_j$. Moreover, both $a_i$ and $a_j$ occur only once.
Hence, applying Lemma~\ref{lemma1} we will not be able to
switch the occurrences of these letters. This and
Proposition~\ref{proposition2} imply that
$e_Xe_Y$ is not an idempotent. This proves the implication
\eqref{proposition70.1}$\Rightarrow$\eqref{proposition70.3}
and completes the proof.
\end{proof}

\begin{corollary}\label{corollary2}
All maximal subgroups of $\mathrm{K}_n$ are trivial
(that is consist of one element).
\end{corollary}

\begin{proof}
Let $f\in \mathrm{K}_n$ be an idempotent and
$x\in \mathrm{K}_n$ be an element, which belongs to the
maximal subgroup of $\mathrm{K}_n$, corresponding to
$f$. Then $x^k=f$ for some $k\in\mathbb{N}$ and
$fx=x^{k+1}=x$. Now Lemma~\ref{lemma3} implies $x=f$,
completing the proof.
\end{proof}

\begin{remark}\label{remark10}
{\rm The idempotent $e_{\{1,\dots,n\}}$ is the
zero element of  $\mathrm{K}_n$. This follows
from Lemma~\ref{lemma1}.
}
\end{remark}

Recall the following {\em natural} order on the idempotents:
$f_1\leq f_2$ if and only if $f_1f_2=f_2f_1=f_1$. We have:

\begin{proposition}\label{propnew1}
Let $f_1,f_2\in \mathrm{K}_n$ be idempotents. Then
$f_1\leq f_2$ if and only if $\mathfrak{c}(f_2)\subset \mathfrak{c}(f_1)$.
\end{proposition}

\begin{proof}
If $\mathfrak{c}(f_2)\subset \mathfrak{c}(f_1)$ then
$f_1f_2=f_2f_1=f_1$ follows from Remark~\ref{remark10}.
Assume that $f_1f_2=f_2f_1=f_1$. Then, by Lemma~\ref{lemma4},
we have $\mathfrak{c}(f_1f_2)=\mathfrak{c}(f_1)\cup\mathfrak{c}(f_2)=
\mathfrak{c}(f_1)$. Hence $\mathfrak{c}(f_2)\subset \mathfrak{c}(f_1)$.
The statement is proved.
\end{proof}

\section{Kiselman's linear representation of $\mathrm{K}_n$}\label{s5}

For $i=1,\dots,n$ denote by $A_i$ the following $(0,1)$-matrix
of size $n\times n$:
\begin{displaymath}
A_i=\left(
\begin{array}{cccccccc}
1&0&\dots &0&1&0&\dots&0\\
0&1&\dots &0&1&0&\dots&0\\
\dots&\dots&\dots &\dots&\dots&\dots&\dots&\dots\\
0&0&\dots &1&1&0&\dots&0\\
0&0&\dots &0&0&0&\dots&0\\
0&0&\dots &0&0&1&\dots&0\\
\dots&\dots&\dots &\dots&\dots&\dots&\dots&\dots\\
0&0&\dots &0&0&0&\dots&1\\
\end{array}
\right),
\end{displaymath}
where the $i$-th row is zero and the $i$-th column 
equals $(1,\dots,1,0,\dots,0)^t$ (the first $i-1$ elements are equal to 
$1$). The following proposition is inspired by \cite[Theorem~3.3]{Ki}.

\begin{proposition}\label{proposition1}
The assignment $a_i\mapsto A_{n-i+1}$ extends uniquely to 
a ho\-mo\-mor\-phism,
$\psi_n:\mathrm{K}_n\to\mathrm{Mat}_{n\times n}(\mathbb{Z})$. Moreover,
we have $\psi_n(e_{\{1,\dots,n\}})=0$.
\end{proposition}

\begin{proof}
Because of \eqref{eq2} it is enough to check that
$A_i^2=A_i$ for all $i=1,\dots,n$; and 
$A_iA_jA_i=A_jA_iA_j=A_jA_i$ for all $i,j$ such that
$1\leq i<j\leq n$. This is a straightforward calculation.
That $\psi_n(e_{\{1,\dots,n\}})=0$ is also a straightforward calculation.
\end{proof}

\begin{remark}\label{remark3}
{\rm
In \cite[Theorem~3.3]{Ki} it is proved that $\psi_3$ is faithful. 
Unfortunately, already $\psi_4$ is not faithful. 
For example, both, $a_3a_4a_2a_1a_3a_2$ 
and $a_3a_2a_4a_3a_1a_2$, are different  canonical words and 
hence represent different elements from  $\mathrm{K}_4$. 
However, one easily computes that $\psi_4(a_3a_4a_2a_1a_3a_2)=
\psi_4(a_3a_2a_4a_3a_1a_2)$.
}
\end{remark}

\section{(Anti)automorphisms of $\mathrm{K}_n$}\label{s6}

\begin{proposition}\label{proposition3}
\begin{enumerate}[(a)]
\item\label{proposition3.1} 
The only automorphism of $\mathrm{K}_n$ is the identity.
\item\label{proposition3.2}
The map $a_i\mapsto a_{n-i+1}$ extends uniquely to an
antiautomorphism of $\mathrm{K}_n$. This is the only
antiautomorphism of $\mathrm{K}_n$.
\end{enumerate}
\end{proposition}

\begin{proof}
Let $\sigma:\mathrm{K}_n\to \mathrm{K}_n$ be an automorphism.
Obviously $\sigma(e)=e$.
The map $\mathfrak{c}\circ \sigma:\mathrm{K}_n\to 2^{\{1,\dots,n\}}$
must be an epimorphism since $\mathfrak{c}$ is an epimorphism 
by Lemma~\ref{lemma4}. For every $i\in \{1,\dots,n\}$
the set $2^{\{1,\dots,n\}}\setminus\{\varnothing,\{i\}\}$ is
closed under $\cup$, and $\mathfrak{c}^{-1}(\{i\})=a_i$. This
implies that $\sigma$ must induce a permutation on
the generators $a_1,\dots,a_n$. Let us prove that 
$\sigma(a_i)=a_i$ by induction on $n$. For $n=1$ the statement is
obvious. By \eqref{eq2}, the letter $a_n$ may be characterized 
as the only letter $a_i$ among $a_1,\dots,a_n$ such that there
does not exist any $a_j$, $j\neq i$, with the property
$a_ja_i=a_ia_ja_i=a_ja_ia_j$. Hence $\sigma(a_n)=a_n$. 
In particular, $\sigma$ induces a permutation of the remaining
letters $a_1,\dots,a_{n-1}$, that is an automorphism of
$\mathrm{K}_{n-1}$. By the inductive assumption, this
automorphism is trivial. Hence $\sigma$ is also trivial.
This proves \eqref{proposition3.1}.

That $a_i\mapsto a_{n-i+1}$ extends uniquely to an
antiautomorphism of $\mathrm{K}_n$ follows from the fact
that it preserves the defining relations \eqref{eq2}.
That this antiautomorphism is unique is proved analogously
to \eqref{proposition3.1}. This completes the proof.
\end{proof}

We will denote the unique antiautomorphism of $\mathrm{K}_n$
by $\tau$.

\begin{question}\label{qnew2}
Is it possible to classify endomorphisms of
$\mathrm{K}_n$?
\end{question}

\section{Green's relations on $\mathrm{K}_n$}\label{s7}

\begin{theorem}\label{theorem4}
Green's relations $\mathcal{L}$, $\mathcal{R}$, 
$\mathcal{D}$, $\mathcal{H}$, and  $\mathcal{J}$ for 
$\mathrm{K}_n$ are trivial (that is all equivalence classes 
of these equivalence relations consist of one element each).
\end{theorem}

To prove this theorem we will need the following notion:
let $A=(a_{i,j})$ be an $n\times n$ matrix with coefficients from some
ring. Define the {\em height} $\mathfrak{h}(A)$ of $A$ as follows:
\begin{displaymath}
\mathfrak{h}(A)=\sum_{i=1}^n |\{j\in\{1,\dots,n\}:a_{i,j}\neq 0\}|\cdot 
2^i.
\end{displaymath}
For $x\in \mathrm{K}_n$ we define the {\em height} $\mathfrak{h}(x)$ 
of $x$ as $\mathfrak{h}(\psi_n(x))$.

We will need the following property of the height:

\begin{lemma}\label{lemma5}
Let $\alpha\in \mathrm{K}_n$ and $i\in \{1,\dots,n\}$
be such that $a_i\alpha\neq \alpha$. Then
$\mathfrak{h}(a_i\alpha)<\mathfrak{h}(\alpha)$. In 
particular, if $\alpha,\beta \in \mathrm{K}_n$ are such that
$\alpha\beta\neq \beta$, then 
$\mathfrak{h}(\alpha\beta)<\mathfrak{h}(\beta)$.
\end{lemma}

\begin{proof}
By the definition of $\mathfrak{h}$ we have to show that
$\mathfrak{h}(\psi_n(a_i)\psi(\alpha))<\mathfrak{h}(\psi_n(\alpha))$.
Set $j=n-i+1$.
Because of the definition of $\psi_n(a_i)=A_{j}$, the matrix
$\psi_n(a_i)\psi_n(\alpha)$ is obtained from the matrix 
$\psi_n(\alpha)$ by the following sequence of elementary
operations: the $j$-th row of $\psi_n(\alpha)$ is added to
all rows with numbers $1,2,\dots,j-1$, and then the 
$j$-th row of the resulting matrix is multiplied with $0$. 
Let $m$ be the number of non-zero entries in the 
$j$-th row of $\psi_n(\alpha)$. This contributes
$m2^{j}$ to $\mathfrak{h}(\alpha)$. Since $\psi_n(\alpha)$
has only non-negative coefficients, adding the $j$-th row
of $\psi_n(\alpha)$ to the rows with numbers $1,2,\dots,j-1$
we can create at most $m$ new non-zero elements in all
these rows. These new elements will contribute at most
$m(2^{j-1}+2^{j-2}+\dots+2^1)<m2^j$ to $\mathfrak{h}(a_i\alpha)$.
Hence $\mathfrak{h}(a_i\alpha)<\mathfrak{h}(\alpha)$ and
the first statement of the lemma is proved. The second
statement follows immediately from the first one.
\end{proof}

Now we are ready to prove Theorem~\ref{theorem4}:

\begin{proof}[Proof of Theorem~\ref{theorem4}.]
Let us prove the statement for the $\mathcal{L}$ relation.
Assume that $a,b\in \mathrm{K}_n$ are such that $a\neq b$ and
$a\mathcal{L}b$. This means that there exists
$x,y\in \mathrm{K}_n$ such that $xa=b$ and $yb=a$.  Hence
from Lemma~\ref{lemma5} we obtain $\mathfrak{h}(b)=
\mathfrak{h}(xa)<\mathfrak{h}(a)$ and
$\mathfrak{h}(a)=\mathfrak{h}(yb)<\mathfrak{h}(b)$.
This implies $\mathfrak{h}(a)<\mathfrak{h}(a)$, a
contradiction. Therefore, every $\mathcal{L}$-class
consists of exactly one element and thus $\mathcal{L}$ 
is trivial.

Since the relation $\mathcal{L}$ is trivial, applying $\tau$ we
obtain that the relation $\mathcal{R}$ is trivial as well.
From the definition of $\mathcal{H}$ and
$\mathcal{D}$ it then follows that both $\mathcal{H}$ and
$\mathcal{D}$ are trivial. Since $\mathrm{K}_n$ is finite,
we have $\mathcal{D}=\mathcal{J}$, completing the proof.
\end{proof}

\begin{remark}\label{remark11}
{\em The statement of Theorem~\ref{theorem4} was announced
in \cite[Theorem~3]{Go}.
}
\end{remark}

\section{Maximal nilpotent subsemigroups of $\mathrm{K}_n$}\label{s8}

Recall that a semigroup, $S$, with the zero element $0$ is called 
{\em nilpotent} provided that there exists $k\in \mathbb{N}$ such that
$S^k=\{0\}$. The minimal possible $k$ with this property is called
the {\em nilpotency class} of $S$. For every $X\subset \{1,\dots,n\}$ 
denote by $\mathrm{Nil}(X)$ the set $\{w\in \mathrm{K}_n|
\mathfrak{c}(w)=X\}$. 

\begin{theorem}\label{theorem5}
\begin{enumerate}[(i)]
\item\label{theorem5.1}
For each $X\subset \{1,\dots,n\}$ the set  $\mathrm{Nil}(X)$ is a 
maximal nilpotent subsemigroup of $\mathrm{K}_n$ 
(with the zero element $e_X$). $\mathrm{Nil}(X)$ has 
nilpotency class $|X|$ if $|X|>0$, and nilpotency class
$1$ if $|X|=0$.
\item\label{theorem5.2}
Every maximal nilpotent subsemigroup of $\mathrm{K}_n$ has the
form $\mathrm{Nil}(X)$ for some $X\subset \{1,\dots,n\}$.
\item\label{theorem5.3}
We have the following decomposition into a disjoint union of maximal 
nilpotent subsemigroups: $\mathrm{K}_n=
\cup_{X\subset \{1,\dots,n\}}\mathrm{Nil}(X)$.
\end{enumerate}
\end{theorem}

\begin{proof}
That $\mathrm{Nil}(X)$ is a subsemigroup of $\mathrm{K}_n$
follows from Lemma~\ref{lemma4}. That $e_X$ is the zero element
of $\mathrm{Nil}(X)$ and the only idempotent of $\mathrm{Nil}(X)$
follows from Lemma~\ref{lemma3}. Hence $\mathrm{Nil}(X)$
is a nilpotent semigroup by \cite[Fact2.30, page~179]{Ar}.
If $w\in \mathrm{K}_n\setminus \mathrm{Nil}(X)$, then 
$w^{|\mathfrak{c}(w)|}$ is an idempotent, different from
$e_X$. This means that the semigroup, generated by $\mathrm{Nil}(X)$
and such $w$, can not be nilpotent. That 
$\mathrm{Nil}(\{\varnothing\})=\{e\}$ has nilpotency class
$1$ is obvious. Let $X\neq \varnothing$. The same arguments as 
the ones used in Lemma~\ref{lemma3} prove that the nilpotency 
class of $\mathrm{Nil}(X)$ is at most $|X|$. Let
$X=\{a_{i_1},\dots,a_{i_k}\}$ and $i_1<i_2<\dots<i_k$. 

\begin{lemma}\label{lemmanilpnew}
The element $w=a_{i_1}a_{i_2}\cdots a_{i_k}$ has order $k$.
\end{lemma}

\begin{proof}
From Lemma~\ref{lemma3} we have that the order of $w$ is at most 
$k$, so we have to prove that $w^l$ is not an idempotent 
for any $l<k$. Observe that, obviously, the subsemigroup of 
$\mathrm{K}_n$, generated by $a_{i_1},a_{i_2},\dots,a_{i_k}$
is isomorphic to $\mathrm{K}_k$ via $a_{i_j}\mapsto a_j$. Hence,
without loss of generality, we may assume
$X=\{1,\dots,n\}$.

By a direct calculation we have that the matrix
$\psi_n(a_1a_2\cdots a_n)$ is an upper triangular matrix
with zero diagonal, whose all element above the diagonal equal $1$.
Hence $\psi_n(a_1a_2\cdots a_n)$ is nilpotent of nilpotency class
exactly $n$. The claim follows.
\end{proof}

From Lemma~\ref{lemmanilpnew} we obtain that
the nilpotency class of $\mathrm{Nil}(X)$ is exactly
$|X|$. This proves \eqref{theorem5.1}.

Let $S$ be a maximal nilpotent subsemigroup of $\mathrm{K}_n$
and $f\in S$ be the corresponding zero element. Then
$f=e_X$ for some $X\subset \{1,\dots,n\}$ by
Proposition~\ref{proposition2}. Since for every element $x$
from $S$ we then should have $x^k=e_X$ for some $k$, from
Lemma~\ref{lemma3} we obtain $S\subset \mathrm{Nil}(X)$,
Now \eqref{theorem5.2} follows from \eqref{theorem5.1}.
The statement \eqref{theorem5.3} is now obvious.
\end{proof}

\section{Isolated and completely isolated subsemigroups of
$\mathrm{K}_n$}\label{s9}

Let $S$ be a semigroup. Recall that a subsemigroup, $T\subset S$,
is called {\em isolated} provided that for all $x\in S$ the
inclusion $x^l\in T$ for some $l\in\mathbb{N}$ implies $x\in T$.
A subsemigroup, $T\subset S$, is called {\em completely isolated} 
provided that $xy\in T$ implies $x\in T$ or $y\in T$ 
for all $x,y\in S$. 

\begin{proposition}\label{proposition21}
\begin{enumerate}[(i)]
\item \label{proposition21.1} The map $\mathfrak{c}$ induces a bijection
between isolated subsemigroups of $\mathrm{K}_n$ and 
subsemigroups of $(2^{\{1,\dots,n\}},\cup)$. In particular, the minimal
isolated subsemigroups of $\mathrm{K}_n$ are $\mathrm{Nil}(X)$,
$X\subset\{1,\dots,n\}$.
\item \label{proposition21.2} The map $\mathfrak{c}$ induces a bijection
between completely isolated subsemigroups of $\mathrm{K}_n$ and
completely isolated  subsemigroups of $(2^{\{1,\dots,n\}},\cup)$.
\end{enumerate}
\end{proposition}

\begin{proof}
Let $S$ be an isolated subsemigroup of $\mathrm{K}_n$. 
Then $\mathfrak{c}(S)=T$ is a subsemigroup of $(2^{\{1,\dots,n\}},\cup)$,
which is obviously isolated since $(2^{\{1,\dots,n\}},\cup)$
consists of idempotents. That $S=\mathfrak{c}^{-1}(T)$
follows from \cite[Proposition~4]{MT}. On the other hand, for any
subsemigroup $T$ of $(2^{\{1,\dots,n\}},\cup)$ the set
$\mathfrak{c}^{-1}(T)$ is a subsemigroup of $\mathrm{K}_n$ and
hence is isolated since $T$ is isolated. This proves 
\eqref{proposition21.1}. \eqref{proposition21.2} follows easily from
\eqref{proposition21.1}.
\end{proof}

\section{Deletion properties}\label{s10}

In this section we establish two combinatorial properties of 
$\mathrm{K}_n$, which  will be used later on during the study 
of linear representations of $\mathrm{K}_n$. However, we think that
these properties are rather remarkable and interesting on their own.

To simplify the notation we set $f=e_{\{2,3,\dots,n\}}$. 
Our {\em first deletion property} is the following statement:

\begin{proposition}\label{proposition15}
Let $v,w\in \mathrm{W}(\{a_2,\dots,a_n\})$ be canonical and different.
Then $va_1f\neq wa_1f$.
\end{proposition}

\begin{proof}
Take the word $va_1f\in \mathrm{W}(\{a_2,\dots,a_n\})$. This word does
not have to be canonical. However, we can use Lemma~\ref{lemma1}
(maybe several times) to reduce it to the unique canonical form
given by Theorem~\ref{theorem2}. Since $v$ is assumed to be canonical,
on the first step we can apply Lemma~\ref{lemma1} only to some subword, 
$a_i\alpha a_i$, of $va_1f$, where the left $a_i$ is a letter of $v$ and
the right $a_i$ is a letter of $f$. This means that $a_1$ is a letter of
$\alpha$, and therefore only Lemma~\ref{lemma1}\eqref{lemma1.1} can
be applied. Thus the new word will have the form $va_1\beta$, where 
$\beta$ is obtained from $f$ by the deletion of one of the letters. The
main point is that the left-hand side $v$ remains the same. Now, 
applying the same argument inductively, we obtain that the canonical
form of $va_1f$ will by $va_1\gamma$, where
$\gamma$ is a quasi-subword of $f$. 

The same argument shows that the canonical form of $wa_1f$ will 
have the form $wa_1\gamma'$, where
$\gamma'$ is a quasi-subword of $f$. Since $a_1$ does not occur in
both $v$ and $w$ by assumption, and $v\neq w$, we obtain that
$va_1\gamma\neq wa_1\gamma'$. The statement now follows from 
Theorem~\ref{theorem2}.
\end{proof}

The {\em second deletion property} is the following more tricky
statement (and is perhaps the deepest result of our paper):

\begin{proposition}\label{proposition16}
Let $w,v,u\in \mathrm{W}(\{a_2,\dots,a_n\})$ be canonical.
Assume that $v\neq u$ and both $wa_1v$ and $wa_1u$ are canonical.
Then $wv\neq wu$, $wva_1\neq wua_1$ and 
$wva_1f\neq wua_1f$.
\end{proposition}

\begin{proof}
We first prove that $wv\neq wu$. Assume this is not the case,
that is assume that $wv=wu$. To proceed we will need some 
preparation.

\begin{lemma}\label{lemma17}
Let $\alpha,\beta\in \mathrm{W}(\{a_2,\dots,a_n\})$ be canonical
and assume that $\alpha a_1\beta$ is canonical as well. Then
the canonical form of $\alpha\beta$ is obtained from
$\alpha\beta$ by deleting some letters of the word $\alpha$
using Lemma~\ref{lemma1}\eqref{lemma1.2}. Moreover, the reduction
process can be organized such that on every step the new
letter which we delete is placed to the left with respect to 
the letter, deleted on the previous step.
\end{lemma}

\begin{proof}
We proceed inductively on the number of deletions. Assume that
$a_i\gamma a_i$ is a subword of $\alpha\beta$, to which we can apply
Lemma~\ref{lemma1}. Since $\alpha a_1\beta$ was canonical,
we obtain that $a_i\gamma a_i=a_i\gamma'\gamma'' a_i$,
where $a_i\gamma'$ is a suffix of $\alpha$ and
$\gamma'' a_i$ is a prefix of $\beta$. Since $a_i\gamma' a_1\gamma'' a_i$,
as a subword of a canonical word, was canonical itself, 
the word $\gamma'\gamma''$ must contain some $a_j$ with $j>i$. Hence
we can only apply Lemma~\ref{lemma1}\eqref{lemma1.2}
to $a_i\gamma a_i$ and thus have to delete some letter from 
$\alpha$. We can of course always start with the rightmost letter of 
$\alpha$, which can be deleted.

Since we delete the rightmost possible letter, the rest of the word,
which is to the right of this letter, has to be canonical. This
part is not affected by our deletion, so it remains canonical. 
On the other  hand, since we have used Lemma~\ref{lemma1}\eqref{lemma1.2}, 
the right neighbor of our letter should have bigger index.
So, if our deletion creates possibilities for new deletions, for these
new possibilities we can only use Lemma~\ref{lemma1}\eqref{lemma1.2}
(this is the same argument as in the previous paragraph).
In particular, it follows that new letters which can be deleted can 
appear only to the left. Moreover, the same argument as above shows that 
if our deletion creates some new letters which can be deleted, it is again 
only Lemma~\ref{lemma1}\eqref{lemma1.2} which can be used. Therefore,
we can again always choose the new rightmost letter and proceed
inductively, completing the proof.
\end{proof}

From Lemma~\ref{lemma17} we obtain that the canonical form
$\mathrm{can}(wu)$ is obtained from $wu$ by deleting some
letters from $w$, and the canonical form
$\mathrm{can}(wv)$ is obtained from $wv$ by deleting some
letters from $w$. In particular, $wu=wv$ implies
$\mathrm{can}(wu)=\mathrm{can}(wv)$. Without loss of generality
we may assume $\mathfrak{l}(u)\leq \mathfrak{l}(v)$. Then the
above observations imply that $v=u'u$ (as a word) for some word $u'$.
In particular, if $\mathfrak{l}(u)= \mathfrak{l}(v)$, we
already get a contradiction, proving that $wv\neq wu$ in this case. 

Hence now we can assume that $\mathfrak{l}(u)< \mathfrak{l}(v)$
and that $v=u'u$ for some non-empty word $u'$. 
Now we are going to make some analysis of
$wu$ and $wv$, which we tried to illustrate on Figure~\ref{fig1}.
It will be convenient for us to distinguish the {\em symbols} 
$\{a_1,\dots,a_n\}$ of our alphabet from the {\em letters} of a
given word (this word will, in fact, be the word $w$). 
So, in the rest of the proof by a {\em letter} of
some word we will mean a symbol of the alphabet together with the
position in the word (so different letters can correspond to
the same symbol). For example, the word $a_1a_2a_3a_1$ is written using
only three different symbols, but it contains four different letters
(the first letter is the symbol $a_1$ staying in position one and the
fourth letter is the the symbol $a_1$ staying in position four).
We will use $a,x,y,v$ to denote the letters
of the words we will work with. 

\begin{figure}
\special{em:linewidth 0.4pt}
\unitlength 0.80mm
\linethickness{0.4pt}
\begin{center}
\begin{picture}(115.00,80.00)
\drawline(7.50,70.00)(65.00,70.00)
\drawline(7.50,60.00)(65.00,60.00)
\drawline(7.50,30.00)(65.00,30.00)
\drawline(7.50,20.00)(65.00,20.00)
\drawline(92.50,70.00)(107.50,70.00)
\drawline(92.50,60.00)(107.50,60.00)
\drawline(92.50,30.00)(107.50,30.00)
\drawline(92.50,20.00)(107.50,20.00)
\drawline(77.50,30.00)(90.00,30.00)
\drawline(77.50,20.00)(90.00,20.00)
\drawline(7.50,70.00)(7.50,60.00)
\drawline(27.50,70.00)(27.50,60.00)
\drawline(32.50,70.00)(32.50,60.00)
\drawline(35.00,70.00)(35.00,60.00)
\drawline(40.00,70.00)(40.00,60.00)
\drawline(42.50,70.00)(42.50,60.00)
\drawline(47.50,70.00)(47.50,60.00)
\drawline(65.00,70.00)(65.00,60.00)
\drawline(92.50,70.00)(92.50,60.00)
\drawline(107.50,70.00)(107.50,60.00)
\drawline(7.50,20.00)(7.50,30.00)
\drawline(27.50,20.00)(27.50,30.00)
\drawline(32.50,20.00)(32.50,30.00)
\drawline(35.00,20.00)(35.00,30.00)
\drawline(40.00,20.00)(40.00,30.00)
\drawline(42.50,20.00)(42.50,30.00)
\drawline(47.50,20.00)(47.50,30.00)
\drawline(50.00,20.00)(50.00,30.00)
\drawline(65.00,20.00)(65.00,30.00)
\drawline(92.50,20.00)(92.50,30.00)
\drawline(107.50,20.00)(107.50,30.00)
\drawline(55.00,20.00)(55.00,30.00)
\drawline(57.50,20.00)(57.50,30.00)
\drawline(62.50,20.00)(62.50,30.00)
\drawline(77.50,20.00)(77.50,30.00)
\drawline(82.50,20.00)(82.50,30.00)
\drawline(90.00,20.00)(90.00,30.00)
\put(80.00,25.00){\makebox(0,0)[cc]{\small{{\it a$'$}}}}
\put(10.00,65.00){\makebox(0,0)[cc]{\small{{\it w}}}}
\put(30.00,65.00){\makebox(0,0)[cc]{\small{{\it x$'$}}}}
\put(37.00,65.00){\makebox(0,0)[cc]{\small{{\it y}}}}
\put(45.00,65.00){\makebox(0,0)[cc]{\small{{\it a}}}}
\put(95.00,65.00){\makebox(0,0)[cc]{\small{{\it u}}}}
\put(10.00,25.00){\makebox(0,0)[cc]{\small{{\it w}}}}
\put(30.00,25.00){\makebox(0,0)[cc]{\small{{\it x$'$}}}}
\put(37.00,25.00){\makebox(0,0)[cc]{\small{{\it y}}}}
\put(45.00,25.00){\makebox(0,0)[cc]{\small{{\it a}}}}
\put(52.50,25.00){\makebox(0,0)[cc]{\small{{\it x}}}}
\put(60.00,25.00){\makebox(0,0)[cc]{\small{{\it y$'$}}}}
\put(87.50,25.00){\makebox(0,0)[cc]{\small{{\it u$'$}}}}
\put(95.00,25.00){\makebox(0,0)[cc]{\small{{\it u}}}}
\put(72.50,42.50){\makebox(0,0)[cc]{$a_1$}}
\dashline{1}(72.50,75.00)(72.50,50.00)
\dashline{1}(72.50,35.00)(72.50,15.00)
\dashline{1}(80.00,30.00)(45.00,60.00)
\dashline{1}(45.00,30.00)(45.00,60.00)
\dashline{1}(52.50,30.00)(30.00,60.00)
\dashline{1}(30.00,30.00)(30.00,60.00)
\dashline{1}(60.0,30.00)(37.50,60.00)
\dashline{1}(37.50,30.00)(37.50,60.00)
\drawline(45.00,60.00)(46.00,58.00)
\drawline(45.00,60.00)(47.20,59.20)
\drawline(45.00,30.00)(46.00,32.00)
\drawline(45.00,30.00)(44.00,32.00)
\drawline(30.00,30.00)(31.00,32.00)
\drawline(30.00,30.00)(29.00,32.00)
\drawline(37.50,60.00)(36.50,58.00)
\drawline(37.50,60.00)(38.50,58.00)
\drawline(30.00,60.00)(30.75,57.70)
\drawline(30.00,60.00)(32.00,58.60)
\drawline(60.00,30.00)(58.00,31.00)
\drawline(60.00,30.00)(59.40,32.20)
\put(56.00,18.00){\makebox(0,0)[cc]{$\underbrace{\hspace{1.3cm}}$}}
\put(56.00,14.00){\makebox(0,0)[cc]{\small{{\it x}}-{\tiny maximal}}}
\put(37.50,18.00){\makebox(0,0)[cc]{$\underbrace{\hspace{0.7cm}}$}}
\put(37.50,14.00){\makebox(0,0)[cc]{\small{{\it y}}-{\tiny maximal}}}
\end{picture}
\end{center}
\caption{Analysis in the proof of Proposition~\ref{proposition16}.}\label{fig1}
\end{figure}
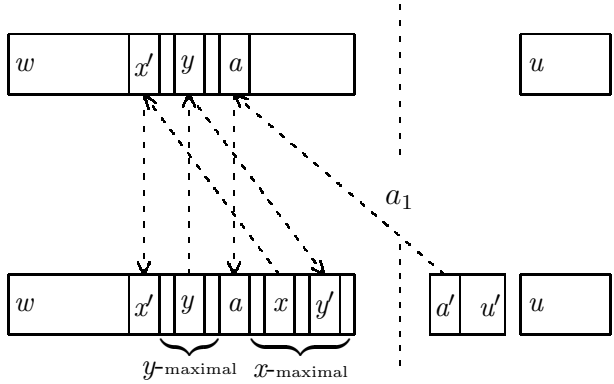

Let $a'$ be the leftmost letter of the non-empty word $u'$. 
Let $a_i$ be the corresponding symbol.
By  Lemma~\ref{lemma17}, the letter $a'$ survives in 
$\mathrm{can}(wv)$. Since $\mathfrak{l}(u)\leq \mathfrak{l}(v)$, 
the corresponding letter of $\mathrm{can}(wu)=\mathrm{can}(wv)$ 
comes from $w$, say from some letter $a$ (this should be one 
of the occurrences of $a_i$ in $w$). Since $wa_1v$ was canonical 
and $a$ is the leftmost letter of $v$,  there should exist a 
symbol, $a_j$, in $w$ to the right of $a$  such that $j>i$. 
We can choose the maximal possible $j$ and let $x$ be the 
rightmost occurrence of $a_j$ in $w$ to the right of our 
letter $a$. All letters in $w$ to the right of $x$ (if any)
have smaller indicies. In $wv$ these letters are followed by
$a'$, which also has smaller index. Hence it is not possible to
delete this $x$ using Lemma~\ref{lemma1}\eqref{lemma1.2}. From
Lemma~\ref{lemma17} we obtain that $x$ survives in 
$\mathrm{can}(wv)$. 

Since $\mathrm{can}(wv)=\mathrm{can}(wu)$, the letter $x$ forces
the existence of some letter $x'$ (representing the same symbol 
$a_j$ as the letter $x$) to the left of $a$, which survives in
$\mathrm{can}(wu)$ and corresponds there to the letter $x$ in
$\mathrm{can}(wv)$. Since $w$ was canonical, between $x'$ and
$x$ in $w$ there should exist some symbol $a_k$ such that $k>j$.
Since $j$ is the maximal possible index to the right of $a$, this
symbol $a_k$ appears in $w$ between $x'$ and $a$. We again
take $k$ the maximal possible and let $y$ be the rightmost
occurrence of $a_k$ between $x'$ and $a$. Then, by definion,
$k$ is bigger than the index of all other symbols
in $w$ to the right of $y$. The letter $a$  survives in
$\mathrm{can}(wu)$, which implies that one can not
use Lemma~\ref{lemma1}\eqref{lemma1.2} to delete $y$ in $wu$.
Hence $y$ survives in $\mathrm{can}(wu)$ between $x'$ and $a$. 

Since $\mathrm{can}(wv)=\mathrm{can}(wu)$, this $y$ should correspond
to some occurrence of $a_k$ to the right of  $x$. However, this
contradicts to the choice of $x$, which was supposed to have the
maximal possible index in $w$ to the right of $a$. The obtained
contradiction proves that $wv=wu$ is not possible, that is the
first inequality of our statement.

Since $wv\neq wu$, the canonical forms $\alpha$ and $\beta$
of $wv$ and $wu$ respectively are different. As $wv,wu\in \mathrm{W}(\{a_2,\dots,a_n\})$ we obtain that  $\alpha a_1$ and
$\beta a_1$ are both canonical and hence different. This proves
the inequality $wva_1\neq wua_1$. From Proposition~\ref{proposition15}
we also obtain $\alpha a_1f\neq \beta a_1f$, which proves
the inequality $wva_1f\neq wua_1f$. This completes the proof.
\end{proof}

\section{Linear representations of $\mathrm{K}_n$}\label{s11}

For a commutative ring, $\mathbf{R}$, we denote by
$\mathbf{R}\mathrm{K}_n$ the semigroup algebra of $\mathrm{K}_n$
over $\mathbf{R}$ and by $\overline{\mathbf{R}\mathrm{K}}_n$ the
quotient of $\mathbf{R}\mathrm{K}_n$ modulo the 
ideal, generated by  the zero element 
$e_{\{1,2,\dots,n\}}$.

\subsection{Faithful representations of $\mathrm{K}_n$}\label{s11.1}

We start with the following observation:

\begin{proposition}\label{proposition22}
Let $\rho$ be a faithful linear representation of $\mathrm{K}_n$
over some field. Then $\dim\rho\geq n$.
\end{proposition}

\begin{proof}
In the proof of Theorem~\ref{theorem5} we saw that the element
$a_1a_2\cdots a_n$ is a nilpotent element of nilpotency class
exactly $n$. Since $e_{\{1,2,\dots,n\}}$ is the zero element in
$\mathrm{K}_n$, factoring, if necessary, the image of
$\rho(e_{\{1,2,\dots,n\}})$ out, we may assume that 
$\rho(e_{\{1,2,\dots,n\}})=0$. If $\rho$ is faithful, the matrix
$\rho(a_1\cdots a_n)$ must then be a nilpotent matrix of
nilpotency class exactly $n$. Obviously, such matrix exists
only if $\dim\rho\geq n$.
\end{proof}

As we have already mentioned in Remark~\ref{remark3}, Kiselman's
representation of  $\mathrm{K}_n$ is not faithful for $n=4$
(and hence for all $n> 4$ either). Let now $\mathbb{K}$ be a
field. From Proposition~\ref{proposition1} we have
$\psi_n(e_{\{1,2,\dots,n\}})=0$ and hence $\psi_n$ is a representation of 
$\overline{\mathbb{K}\mathrm{K}}_n$ as well. We continue with the following
observation about faithfulness:

\begin{proposition}\label{proposition23}
The indecomposable projective cover of Kiselman's representation 
of $\overline{\mathbb{K}\mathrm{K}}_n$ in $\mathbb{K}^n$ is faithful
as a representation of $\mathrm{K}_n$.
\end{proposition}

\begin{proof}
Set $\pi_1=e-a_n\in \mathbb{K}\mathrm{K}_n$,
$\pi_2=a_n-a_na_{n-1}\in \mathbb{K}\mathrm{K}_n$,\dots,
$\pi_{n-1}=a_na_{n-1}\cdots a_{3}-a_na_{n-1}\cdots a_2
\in \mathbb{K}\mathrm{K}_n$, $\pi_{n}=a_na_{n-1}\cdots a_2$.
By a direct calculation using the formulae from Section~\ref{s5}
one obtains that for $i=1,\dots,n$ the matrix
$\psi_n(\pi_i)$ is the diagonal matrix $D_i$, whose diagonal is
the vector  $(0,\dots,0,1,0,\dots,0)$, where the element $1$ stays on the 
$i$-th place.

First we claim that the vector $v=(0,0,\dots,0,1)^t$ generates
Kiselman's representation. Indeed, $A_1v=(1,1,\dots,1,0)^t$ and hence,
acting on $A_1v$ by $D_i$, $i=1,\dots,n-1$, we produce 
all elements from the standard basis of $\mathbb{K}^n$. 

From Proposition~\ref{proposition2} we know that $\pi_n=
e_{\{2,3,\dots,n\}}$ is
an idempotent. Furthermore, $\psi_n(\pi_n)v=v$ and hence 
$\overline{\mathbb{K}\mathrm{K}}_n\pi_n$ is a projective cover of
Kiselman's representation.

Every element of $\mathrm{K}_n$ can be written as either $w$ or
$wa_1v$, where $w,v\in \mathrm{W}(\{a_2,\dots,a_n\})$. From 
Remark~\ref{remark10} it follows that $\pi_n w=w\pi_n=
\pi_n v=v\pi_n=\pi_n$. Hence for any $\alpha\in \mathrm{K}_n$
we have
\begin{displaymath}
\pi_n\alpha\pi_n=
\begin{cases}
\pi_n,& \text{ $a_1$ is not a letter of $\alpha$};\\
e_{\{1,2,\dots,n\}},& \text{ otherwise}.
\end{cases}
\end{displaymath}
Hence $\pi_n\mathbb{K}\mathrm{K}_n\pi_n$ has dimension two and
a monomial basis, consisting of $\pi_n$ and
$e_{\{1,2,\dots,n\}}$. Factoring out the zero element
$e_{\{1,2,\dots,n\}}$ we get a copy of the ground field since
$\pi_n$ is an idempotent. Thus 
$\pi_n\overline{\mathbb{K}\mathrm{K}}_n\pi_n$ is a local algebra. 
Hence $\pi_n$ is a primitive idempotent of 
$\overline{\mathbb{K}\mathrm{K}}_n$,
which implies that the $\overline{\mathbb{K}\mathrm{K}}_n$-module
$\overline{\mathbb{K}\mathrm{K}}_n\pi_n$ is indecomposable.

To complete the proof we have just to show that the corresponding
representation of $\mathrm{K}_n$ is faithful. By definition,
the module $\overline{\mathbb{K}\mathrm{K}}_n\pi_n$ has a 
monomial basis, which consists of all non-zero elements from the 
left principal ideal of $\mathrm{K}_n$, generated by $\pi_n$. In 
particular, we have the basis elements $\pi_n$ and $a_1\pi_n$ 
(note that $a_1\pi_n$ is a canonical word).

If $w,v\in \mathrm{W}(\{a_2,\dots,a_n\})$ are different and canonical, 
then $wa_1\pi_n\neq va_1\pi_n$ by Proposition~\ref{proposition15}.
The elements $wa_1\pi_n$ and $va_1\pi_n$ are linearly independent in
$\mathbb{K}\mathrm{K}_n\pi_n$, in particular, they are different.
Therefore the elements $w$ and $v$ from $\mathrm{K}_n$ are
represented by different linear operators on 
$\mathbb{K}\mathrm{K}_n\pi_n$.

If $u,v,w\in \mathrm{W}(\{a_2,\dots,a_n\})$ are canonical, then
$u\pi_n=\pi_n$ and $va_1w\pi_n=va_1\pi_n\neq \pi_n$. Hence
the elements $u$ and $va_1w$ from $\mathrm{K}_n$ are
represented by different linear operators on 
$\overline{\mathbb{K}\mathrm{K}}_n\pi_n$.

Let $w_1a_1v_1$ and $w_2a_1v_2$ be two different elements from
$\mathrm{K}_n$, written in the canonical form. In particular,
$w_1,w_2,v_1,v_2\in \mathrm{W}(\{a_2,\dots,a_n\})$ and are canonical.
If $w_1\neq w_2$, we have $w_1a_1v_1\pi_n=w_1a_1\pi_n$ and
$w_2a_1v_2\pi_n=w_2a_1\pi_n$ (since $\pi_n$ is the zero element with
respect to $a_j$, $j>1$). Moreover, from Proposition~\ref{proposition15}
we get $w_1a_1\pi_n\neq w_2a_1\pi_n$. Both $w_1a_1\pi_n$ and
$w_2a_1\pi_n$ are basis elements of $\overline{\mathbb{K}\mathrm{K}}_n\pi_n$,
which implies that the elements $w_1a_1v_1$ and $w_2a_1v_2$ are
represented by different linear operators on
$\overline{\mathbb{K}\mathrm{K}}_n\pi_n$.

Assume now that $w_1=w_2=w$. Then $v_1\neq v_2$ and we have
$w_1a_1v_1a_1\pi_n=w_1v_1a_1\pi_n$ and
$w_2a_1v_2a_1\pi_n=w_2v_2a_1\pi_n$ using 
Lemma~\ref{lemma1}\eqref{lemma1.2}. From Proposition~\ref{proposition16}
we get $w_1v_1a_1\pi_n\neq w_2v_2a_1\pi_n$. Both $w_1v_1a_1\pi_n$ and
$w_2v_2a_1\pi_n$ are basis elements of 
$\overline{\mathbb{K}\mathrm{K}}_n\pi_n$,
which implies that the elements $w_1a_1v_1$ and $w_2a_1v_2$ are
represented by different linear operators on 
$\overline{\mathbb{K}\mathrm{K}}_n\pi_n$.
Hence the representation of $\mathrm{K}_n$ on
$\overline{\mathbb{K}\mathrm{K}}_n\pi_n$ is faithful.
\end{proof}

The ideas from the proof of Proposition~\ref{proposition23} can
be used to construct a huge family of faithful $n$-dimensional
representations of $\mathrm{K}_n$. Consider the 
polynomial ring $\mathbb{Z}[\xi_{i,j}:1\leq i<j\leq n]$. 
Define the following representation of $\mathrm{K}_n$ by 
$n\times n$-matrices over $\mathbb{Z}[\xi_{i,j}:1\leq i<j\leq n]$:
\begin{displaymath}
\kappa_n:a_{n-i+1}\mapsto \left(
\begin{array}{cccccccc}
1&0&\dots &0&\xi_{1,i}&0&\dots&0\\
0&1&\dots &0&\xi_{2,i}&0&\dots&0\\
\dots&\dots&\dots &\dots&\dots&\dots&\dots&\dots\\
0&0&\dots &1&\xi_{i-1,i}&0&\dots&0\\
0&0&\dots &0&0&0&\dots&0\\
0&0&\dots &0&0&1&\dots&0\\
\dots&\dots&\dots &\dots&\dots&\dots&\dots&\dots\\
0&0&\dots &0&0&0&\dots&1\\
\end{array}
\right),
\end{displaymath}
where the $i$-th row is zero and the $i$-th column 
equals $(\xi_{1,i},\dots,\xi_{i-1,i},0,\dots,0)^t$.

\begin{proposition}\label{proposition25}
The representation $\kappa_n$ is faithful.
\end{proposition}

\begin{proof}
We proceed by induction on $n$. For $n=1,2$ the statement is easily
checked by a direct calculation. 

Let $w,u\in \mathrm{W}(\{a_2,\dots,a_n\})$ be different and canonical.
The semigroup generated by $a_2,\dots,a_n$ is obviously isomorphic
to $\mathrm{K}_{n-1}$ under the map $a_i\mapsto a_{i-1}$. Let us
denote this isomorphism by $F$. Then the first $n-1$ rows and the 
first $n-1$ columns of  $\kappa_n(w)$ and $\kappa_n(u)$ are 
exactly the matrices $\kappa_{n-1}(F(w))$ and 
$\kappa_{n-1}(F(u))$ respectively. By induction we have
$\kappa_{n-1}(F(w))\neq \kappa_{n-1}(F(u))$ and hence
$\kappa_n(w)\neq\kappa_n(u)$.

Let $u,v,w\in \mathrm{W}(\{a_2,\dots,a_n\})$ be canonical. Then
the last diagonal element of $u$ is $1$ while the last diagonal
element of $va_1w$ is $0$. Hence $\kappa_n(u)\neq\kappa_n(va_1w)$.

Let $w_1,w_2,v_1,v_2\in \mathrm{W}(\{a_2,\dots,a_n\})$ be canonical.
Assume that $w_1\neq w_2$ and that $w_1a_1v_1$ and $w_2a_1v_2$
are also canonical. Recall that $\pi_n=a_n\cdots a_2$.
As in the proof of Proposition~\ref{proposition23}
we have $w_1a_1v_1\pi_n=w_1a_1\pi_n$ and $w_2a_1v_2\pi_n=w_2a_1\pi_n$.
Further
\begin{displaymath}
\kappa_n(a_1\pi_n)=:\left(
\begin{array}{ccccc}
0&0&\dots&0&\xi_{1,n}\\
0&0&\dots&0&\xi_{2,n}\\
\dots&\dots&\dots&\dots&\dots\\
0&0&\dots&0&\xi_{n-1,n}\\
0&0&\dots&0&0
\end{array}
\right).
\end{displaymath}
Since $w_1\neq w_2$, by induction we, similarly to the arguments above,
derive that the matrices $M_1$ and $M_2$, formed by the first $n-1$ 
rows and the  first $n-1$ columns of the matrices $\kappa_n(w_1)$ 
and $\kappa_n(w_2)$ respectively, are different. Since
$w_1,w_2\in \mathrm{W}(\{a_2,\dots,a_n\})$, the coefficients of these
matrices do not contain $\xi_{i,n}$ for all $i$. Now observe that
$\xi_{1,n},\dots,\xi_{n-1,n}$ are linearly independent 
(over $\mathbb{Z}[\xi_{i,j}:1\leq i<j\leq n-1]$) elements
of $\mathbf{R}$. From the definition of the matrix multiplication we
get that the last columns in the matrices $\kappa_n(w_1a_1\pi_n)$
and $\kappa_n(w_2a_1\pi_n)$ will be different. Hence
$\kappa_n(w_1a_1v_1\pi_n)\neq\kappa_n(w_2a_1v_1\pi_n)$ and therefore
$\kappa_n(w_1a_1v_1)\neq\kappa_n(w_2a_1v_1)$.

Finally, let us assume that 
$w,u,v\in \mathrm{W}(\{a_2,\dots,a_n\})$ are canonical and
such that $wa_1u$ and $wa_1v$ are canonical and different.
By Lemma~\ref{lemma1}\eqref{lemma1.2} we have 
$wa_1ua_1\pi_n=wua_1\pi_n$ and $wa_1va_1\pi_n=wva_1\pi_n$.
Moreover, from Proposition~\ref{proposition25} we have
$wu\neq wv$. The same arguments as in the previous paragraph show that
the last columns in the matrices $\kappa_n(wua_1\pi_n)$
and $\kappa_n(wua_1\pi_n)$ will be different. Hence
$\kappa_n(wa_1u)\neq\kappa_n(wa_1v)$. This completes the proof.
\end{proof}

As an immediate corollary we obtain the following statement, which,
together with Proposition~\ref{proposition22}, was announced in 
\cite[Theorem~4]{Go}:

\begin{theorem}\label{theorem27}
$\mathrm{K}_n$ has a faithful representation by 
$n\times n$ matrices with non-negative integer coefficients.
\end{theorem}

\begin{proof}
By Proposition~\ref{proposition25}, the representation 
$\kappa_n$ is faithful. For every pair
$\{\alpha,\beta\}$ of different elements from
$\mathrm{K}_n$ we have $\kappa_n(\alpha)\neq \kappa_n(\beta)$,
hence there exist $i_{\{\alpha,\beta\}}$ and $j_{\{\alpha,\beta\}}$
such that the $(i_{\{\alpha,\beta\}},j_{\{\alpha,\beta\}})$-entry
of $\kappa_n(\alpha)$ is different from the
$(i_{\{\alpha,\beta\}},j_{\{\alpha,\beta\}})$-entry 
of $\kappa_n(\beta)$. These entries are polynomials with integer
coefficients, so this condition can be written as the condition
``some non-zero polynomial in $\xi_{i,j}$ is not equal to zero''.
Since $\mathrm{K}_n$ is finite by Theorem~\ref{theorem1}, the
faithfullness of $\kappa_n$ gives us a finite number of
polynomial inequalities. Since the set
$\mathbb{N}^{n(n-1)/2}$ is Zariski dense in
$\mathbb{Q}^{n(n-1)/2}$, we will get that there are infinitely
many collections of $n_{i,j}\in \mathbb{N}$, $1\leq i<j\leq n$,
such that after the evaluation $\xi_{i,j}\to n_{i,j}$ all
our inequalities are still satisfied. This means that there 
are infinitely many collections of $n_{i,j}\in \mathbb{N}$, 
$1\leq i<j\leq n$, such that after the evaluation 
$\xi_{i,j}\to n_{i,j}$ we obtain a faithful representation of
$\mathrm{K}_n$ with non-negative integer coefficients.
This completes the proof.
\end{proof}

Following the proof of Proposition~\ref{proposition25} one
can in fact explicitly present a collection of $n_{i,j}$,
such that after the evaluation $\xi_{i,j}\to n_{i,j}$
one obtains a faithful representation of
$\mathrm{K}_n$ with non-negative integer coefficients.
Define two sequences, $\mathrm{m}_i$ and $\mathrm{l}_i$,
$i\geq 1$, recursively as follows: $\mathrm{m}_1=\mathrm{l}_1=1$,
$\mathrm{m}_i=\mathrm{l}_{i-1}+1$, $\mathrm{l}_i=i^{2^i}
\mathrm{m}_{i}^{i2^i}$, $i\geq 2$.

\begin{proposition}\label{proposition31}
Denote by $\kappa'_n$ the representation of $\mathrm{K}_n$ with 
non-negative integer coefficients, obtained from $\kappa_n$
via the evaluation $\xi_{i,j}\to \mathrm{m}_j^i$.
\begin{enumerate}[(i)]
\item\label{proposition31.1} $\kappa'_n$ is faithfull.
\item\label{proposition31.2} For every $w\in \mathrm{K}_n$ 
each entry of the matrix $\kappa'_n(w)$ is smaller than $\mathrm{l}_n$.
\end{enumerate}
\end{proposition}

\begin{proof}
We prove this by the simultaneous induction on $n$. For
$n=2$ both statements are easily checked by a direct calculation.
Since $\mathrm{m}_{n}^i>\mathrm{l}_j$ for all $i\geq 1$ and
$j<n$ by construction, the maximal possible entry appearing
in the matrix $\kappa_{n}'(a_i)$, $i\leq n$, is
$\mathrm{m}_{n}^{n-1}<\mathrm{m}_{n}^n$. From 
Corollary~\ref{cnew1}\eqref{cnew1.3} it follows that every element
from $\mathrm{K}_{n}$ can be written as a product of at
most $2^{n}$ generators. It is easy to see that then
the maximal possible entry of such product is smaller than
$n^{2^n}(\mathrm{m}_{n}^{n})^{2^n}$. The induction
step for \eqref{proposition31.2} is now completed by comparing 
this with the definition of $\mathrm{l}_n$. 

To prove \eqref{proposition31.1} we just follow the proof of 
Proposition~\ref{proposition25}. It is easy to see that the
only thing we have to verify is that, given two different matrices
$\kappa_{n-1}'(F(w))$ and $\kappa_{n-1}'(F(v))$, the rightmost columns
of the matrices $\kappa_{n}(wa_1\pi_n)$ and
$\kappa_{n}(va_1\pi_n)$ are different. These columns are linear
combinations of $\mathrm{m}_n^i$, $i=1,\dots,n-1$ with coefficients
from the matrices $\kappa_{n-1}'(F(w))$ and $\kappa_{n-1}'(F(u))$.
By induction, all such coefficients do not exceed
$\mathrm{l}_{n-1}$, which is strictly smaller than $\mathrm{m}_n$
by definition. It follows that two such linear combinations with 
different collections of such coefficients will be different. This
completes the proof.
\end{proof}

\subsection{Irreducible representations and the structure
of $\mathbb{K}\mathrm{K}_n$}\label{s11.2}

Let $\mathbb{K}$ be a field. For any $X\subset \{1,2,\dots,n\}$ we
define the map $\rho_{X}:\mathrm{K}_n\to \mathbb{K}$ as follows:
\begin{displaymath}
\rho_X(w)=
\begin{cases}
1,& \mathfrak{c}(w)\subset X;\\
0,& \text{otherwise}.
\end{cases}
\end{displaymath}

\begin{proposition}\label{proposition41}
\begin{enumerate}[(i)]
\item\label{proposition41.1} For any $X\subset \{1,2,\dots,n\}$ 
the map $\rho_X$ gives an irreducible 
representation of $\mathbb{K}\mathrm{K}_n$.
\item\label{proposition41.2} Representations $\rho_X$, 
$X\subset \{1,2,\dots,n\}$, are pairwise non-equivalent and constitute an
exhaustive list of irreducible representations of $\mathbb{K}\mathrm{K}_n$.
In particular, $\mathbb{K}\mathrm{K}_n$ has $2^n$ non-equivalent 
irreducible representations.
\item\label{proposition41.3} $\rho_X$ is a representation of
$\overline{\mathbb{K}\mathrm{K}}_n$ if and only if 
$X\neq \{1,2,\dots,n\}$. In particular, $\overline{\mathbb{K}\mathrm{K}}_n$
has $2^{n}-1$ non-equivalent  irreducible representations.
\end{enumerate}
\end{proposition}

\begin{proof}
Fix $X\subset \{1,2,\dots,n\}$. For $i\in\{1,\dots,n\}$
define $\overline{\rho}_X(a_i)$ to be $1$ if $i\in X$ and $0$ otherwise.
It is straightforward to check that this assignment satisfies the defining
relations \eqref{eq2} of $\mathrm{K}_n$. Hence it extends uniquely 
to a representation of $\mathrm{K}_n$. From the definition of
$\mathfrak{c}$ one immediately obtains that this extension is the
map $\rho_X$. The representation $\rho_X$ is irreducible since it is
one-dimensional. This proves \eqref{proposition41.1}.

Let $X$ and $Y$ be different subsets of $\{1,\dots,n\}$. Withour loss
of generality we may assume that $X\setminus Y\neq\varnothing$. Let
$i\in X\setminus Y$. Then $\rho_X(a_i)=1$ and $\rho_Y(a_i)=0$. Hence
$\rho_X$ and $\rho_Y$ are not equivalent. In particular, we have
$2^n$ non-equivalent irreducible representations of 
$\mathbb{K}\mathrm{K}_n$. However, from Proposition~\ref{proposition2}
we know that $\mathrm{K}_n$ has $2^n$ idempotents, and from 
Theorem~\ref{theorem4} we know that all Green's relations on
$\mathrm{K}_n$ are trivial. Hence, Munn's Theorem (see for example
\cite[Theorem~5.33]{CP}) gives us that $\mathbb{K}\mathrm{K}_n$
has exactly $2^n$ non-equivalent irreducible representations.
This proves \eqref{proposition41.2}. \eqref{proposition41.3}
follows immediately from \eqref{proposition41.1},
\eqref{proposition41.2} and a direct calculation. This completes the
proof.
\end{proof}

\begin{corollary}\label{corollary44}
The algebra $\mathbb{K}\mathrm{K}_n$ is basic.
\end{corollary}

\begin{proof}
From Proposition~\ref{proposition41}\eqref{proposition41.2}
we have that all simple $\mathbb{K}\mathrm{K}_n$-modules are
one-dimensional. This implies the statement.
\end{proof}

Since we now know all irreducible representations of 
$\mathbb{K}\mathrm{K}_n$, it is a natural question to determine 
the decomposition of the regular module into a direct sum of
indecomposable projectives, that is to find a decomposition of the
unit element of $\mathbb{K}\mathrm{K}_n$ into a direct sum of 
pairwise orthogonal primitive idempotents.

Let $X\subset \{1,\dots,n\}$. Assume that 
$X=\{i_1,\dots,i_s\}$, where $i_1>i_2>\dots>i_s$;
and $\{1,\dots,n\}\setminus X=\{j_1,\dots,j_t\}$, where
$j_1<j_2<\dots<j_t$. Set 
\begin{displaymath}
e_{X}^{(n)}=a_{i_1}a_{i_2}\cdots a_{i_s}
(e-a_{j_1})(e-a_{j_2})\cdots(e-a_{j_t})\in \mathbb{K}\mathrm{K}_n.
\end{displaymath}

\begin{proposition}\label{proposition42}
\begin{enumerate}[(i)]
\item\label{proposition42.1}
\begin{multline*}
\{e_{X}^{(n)}\,:\, X \subset \{1,\dots,n\}\}=
a_n\{e_{Y}^{(n-1)}\,:\, Y \subset \{1,\dots,n-1\}\}\cup\\ \cup
 \{e_{Y}^{(n-1)}\,:\, Y \subset \{1,\dots,n-1\}\}(e-a_n).
\end{multline*}
\item\label{proposition42.2} 
For every $X\subset \{1,\dots,n\}$ the element
$e_{X}^{(n)}$ is a primitive idempotent of $\mathbb{K}\mathrm{K}_n$.
\item\label{proposition42.3} 
$e_{X}^{(n)}e_{Y}^{(n)}=0$ if $X\neq Y$.
\item\label{proposition42.4}
$e=\sum_{X\subset \{1,\dots,n\}}e_{X}^{(n)}$. 
\end{enumerate}
\end{proposition}

\begin{proof}
If $n\in X$, from the definition of  $e_{X}^{(n)}$ we have
$e_{X}^{(n)}=a_ne_{X\setminus\{n\}}^{(n-1)}$. 
If $n\not\in X$, from the definition of  $e_{X}^{(n)}$ we have
$e_{X}^{(n)}=e_{X}^{(n-1)}(e-a_n)$. This proves \eqref{proposition42.1}. 

Now we prove the rest by a simultaneous induction on $n$.
For $n=1$ the statements \eqref{proposition42.2}, 
\eqref{proposition42.3} and \eqref{proposition42.4} are obvious.

Let $Y \subset \{1,\dots,n-1\}$. Then
\begin{displaymath}
\begin{array}{lcl}
a_ne_{Y}^{(n-1)}a_ne_{Y}^{(n-1)}&=&
\text{ (by Lemma~\ref{lemma1}\eqref{lemma1.1}) }\\
a_ne_{Y}^{(n-1)}e_{Y}^{(n-1)}&=&
\text{ (by inductive assumption) }\\
a_ne_{Y}^{(n-1)}.&&
\end{array}
\end{displaymath}
Analogously, using Lemma~\ref{lemma1}\eqref{lemma1.1} and
the inductive assumption, we have
\begin{displaymath}
\begin{array}{lc}
e_{Y}^{(n-1)}(e-a_n)e_{Y}^{(n-1)}(e-a_n)&=\\
e_{Y}^{(n-1)}e_{Y}^{(n-1)}
-e_{Y}^{(n-1)}a_ne_{Y}^{(n-1)}-
e_{Y}^{(n-1)}e_{Y}^{(n-1)}a_n+
e_{Y}^{(n-1)}a_ne_{Y}^{(n-1)}a_n&=\\
e_{Y}^{(n-1)}e_{Y}^{(n-1)}
-e_{Y}^{(n-1)}a_ne_{Y}^{(n-1)}-
e_{Y}^{(n-1)}e_{Y}^{(n-1)}a_n+
e_{Y}^{(n-1)}a_ne_{Y}^{(n-1)}&=\\
e_{Y}^{(n-1)}- e_{Y}^{(n-1)}a_n&=\\
e_{Y}^{(n-1)}(e-a_n).&
\end{array}
\end{displaymath}
Hence all $e_{X}^{(n)}$ are idempotents. 

Let $Y,Z\subset \{1,\dots,n-1\}$. Then, using
Lemma~\ref{lemma1}\eqref{lemma1.1} and the inductive 
assumption, we compute:
\begin{gather*}
a_ne_{Y}^{(n-1)}a_ne_{Z}^{(n-1)}=a_ne_{Y}^{(n-1)}e_{Z}^{(n-1)}=0;\\
a_ne_{Y}^{(n-1)}e_{Z}^{(n-1)}(e-a_n)=0;\\
e_{Z}^{(n-1)}(e-a_n)a_ne_{Y}^{(n-1)}=0.
\end{gather*}
Finally,
\begin{multline*}
e_{Y}^{(n-1)}(e-a_n)e_{Z}^{(n-1)}(e-a_n)=\\=
e_{Y}^{(n-1)}e_{Z}^{(n-1)}-
e_{Y}^{(n-1)}a_ne_{Z}^{(n-1)}-
e_{Y}^{(n-1)}e_{Z}^{(n-1)}a_n+
e_{Y}^{(n-1)}a_ne_{Z}^{(n-1)}a_n=\\
=-e_{Y}^{(n-1)}a_ne_{Z}^{(n-1)}+e_{Y}^{(n-1)}a_ne_{Z}^{(n-1)}=0.
\end{multline*}
Hence the idempotents $e_{X}^{(n)}$, $X\subset \{1,\dots,n\}$,
are pairwise orthogonal.

Further, using \eqref{proposition42.1} and the inductive assumption
we have
\begin{multline*}
\sum_{X\subset \{1,\dots,n\}}e_{X}^{(n)}=
a_n\left(\sum_{Y\subset \{1,\dots,n-1\}}e_{Y}^{(n-1)}\right)+
\left(\sum_{Y\subset \{1,\dots,n-1\}}e_{Y}^{(n-1)}\right)(e-a_n)=\\
= a_n+(e-a_n)=e.
\end{multline*}

By the definition of $e_{X}^{(n)}$, the element
$e_{X}^{(n)}$ is a linear combination of different canonical
monomials. Hence $e_{X}^{(n)}\neq 0$ in $\mathbb{K}\mathrm{K}_n$.
Now since the number of different $e_{X}^{(n)}$'s is $2^n$,
the statement about the primitivity of $e_{X}^{(n)}$'s
follows from Proposition~\ref{proposition41}\eqref{proposition41.2}
and Corollary~\ref{corollary44}. This completes the proof.
\end{proof}

\begin{corollary}\label{corollary46}
Let $X\subset\{1,2,\dots,n\}$. Then
$\mathbb{K}\mathrm{K}_ne_{X}^{(n)}$ is the projective cover
of $\rho_X$.
\end{corollary}

\begin{proof}
It is a straightforward calculation that 
$\rho_X(e_{X}^{(n)})=1$. The claim follows.
\end{proof}

\begin{remark}\label{remark47}
{\rm
One easily checks that the simple subquotients of 
Kiselman's representation of $\mathbb{K}\mathrm{K}_n$
are $\rho_X$, where $|\{1,2,\dots,n\}\setminus X|=1$,
each occurring with multiplicity one.
}
\end{remark}

As one more immediate corollary we obtain the following very
surprising result, which once more emphasizes the importance of
Kiselman's representation and shows that 
Proposition~\ref{proposition23} is fairly remarkable:

\begin{corollary}\label{corollary48}
Let $X\subset \{1,2,\dots,n\}$ be such that
$X\neq \{2,3,\dots,n\}$. Then the projective module $\mathbb{K}\mathrm{K}_ne_{X}^{(n)}$
is {\em not} a faithful representation of
$\mathrm{K}_n$.
\end{corollary}

\begin{proof}
The statement is obvious in the case $X=\{1,2,\dots,n\}$,
so we may assume $X\neq \{1,2,\dots,n\}$.
Set $w=e_{\{2,3,\dots,n\}}-e_{\{1,2,\dots,n\}}\in 
\mathbb{K}\mathrm{K}_n$. It is certainly enough to show that 
$w\mathbb{K}\mathrm{K}_ne_{X}^{(n)}=0$
(which means that the different elements $e_{\{2,3,\dots,n\}}$
and $e_{\{1,2,\dots,n\}}$ are represented by the same
linear transformations on $\mathbb{K}\mathrm{K}_ne_{X}^{(n)}$).
For $v\in \mathrm{W}(\{a_1,\dots,a_n\})$ we have
\begin{displaymath}
wv=
\begin{cases}
w,&\text{ $v$ does not contain $a_1$};\\
e_{\{1,2,\dots,n\}}, & \text{otherwise}.
\end{cases}
\end{displaymath}
Hence for any $x\in \mathbb{K}\mathrm{K}_n$
we have $wx=\alpha w+\beta e_{\{1,2,\dots,n\}}$ for some
$\alpha,\beta\in\mathbb{K}$. Therefore
\begin{displaymath}
wxe_{X}^{(n)}=\alpha we_{X}^{(n)}+\beta e_{\{1,2,\dots,n\}}e_{X}^{(n)}
=
\alpha e_{\{2,3,\dots,n\}}^{(n)}e_{X}^{(n)}+
\beta e_{\{1,2,\dots,n\}}^{(n)}e_{X}^{(n)}=0
\end{displaymath}
by Proposition~\ref{proposition42}\eqref{proposition42.3}. The claim
follows.
\end{proof}

One can now say even more about the structure of $\mathbb{K}\mathrm{K}_n$,
in particular, giving an independent explanation for 
Corollary~\ref{corollary48}:

\begin{proposition}\label{proposition49}
The algebra $\mathbb{K}\mathrm{K}_n$ is directed in the sense that there
exists a linear order, $\prec$, on the set
$\{X:X\subset \{1,2,\dots,n\}\}$ such that
\begin{displaymath}
\mathrm{Hom}_{\mathbb{K}\mathrm{K}_n}
(\mathbb{K}\mathrm{K}_ne_{X}^{(n)},\mathbb{K}\mathrm{K}_ne_{Y}^{(n)})=0
\end{displaymath}
provided that $Y\prec X$. In particular, the algebra 
$\mathbb{K}\mathrm{K}_n$ is quasi-hereditary with respect to
$\prec$ with projective standard modules.
\end{proposition}

\begin{proof}
Let us prove directness by induction on $n$. For $n=1$ the statement 
is obvious. To prove the induction step we consider the projective 
modules $P_1=\mathbb{K}\mathrm{K}_n a_n$ and $P_2=
\mathbb{K}\mathrm{K}_n (e-a_n)$. Obviously
$\mathbb{K}\mathrm{K}_n\cong P_1\oplus P_2$.

Observe that for any $x\in \mathrm{K}_n$, using
Lemma~\ref{lemma1}\eqref{lemma1.1}, we have
\begin{displaymath}
a_nx(e-a_n)=a_nx-a_nxa_n=a_nx-a_nx=0.
\end{displaymath}
Hence $\mathrm{Hom}_{\mathbb{K}\mathrm{K}_n}(P_1,P_2)=0$.

The endomorphism algebra of $P_1$ is the opposite of the algebra 
$B=a_n\mathbb{K}\mathrm{K}_n a_n$. This algebra is the linear span of
the set $\{a_nxa_n:x\in \mathrm{K}_n\}$. Using
Lemma~\ref{lemma1}\eqref{lemma1.1}, every element from the latter set can
be written  as $a_ny$, where $y\in \mathrm{K}_{n-1}$, moreover all
such elements are obviously linearly independent. It follows that
$a_ny\mapsto y$ induces an isomorphism of $B$ onto 
$\mathbb{K}\mathrm{K}_{n-1}$.
By the inductive assumption we obtain that $B$ is directed.

The endomorphism algebra of $P_2$ is the opposite of the algebra 
$C=(e-a_n)\mathbb{K}\mathrm{K}_n (e-a_n)$. This algebra is the linear span of
the set $\{(e-a_n)x(e-a_n):x\in \mathrm{K}_n\}$. Note that
\begin{displaymath}
(e-a_n)x(e-a_n)=x-a_nx-xa_n+a_nxa_n=x-xa_n
\end{displaymath}
by Lemma~\ref{lemma1}\eqref{lemma1.1}. In particular, if
$x$ contains $a_n$, then from Lemma~\ref{lemma1}\eqref{lemma1.1}
it follows that $(e-a_n)x(e-a_n)=x-xa_n=x-x=0$.
This means that $C$ has the following basis:
$\{(e-a_n)x(e-a_n):x\in \mathrm{K}_{n-1}\}$ and one immediately
checks that $(e-a_n)x(e-a_n)\mapsto x$ induces an isomorphism
from $C$ onto $\mathbb{K}\mathrm{K}_{n-1}$. By the inductive assumption 
we obtain that $C$ is directed as well.

So, the endomorphism algebras of both $P_1$ and $P_2$ are directed
and $\mathrm{Hom}_{\mathbb{K}\mathrm{K}_n}(P_1,P_2)=0$. It follows
that $\mathbb{K}\mathrm{K}_n$ is directed, as asserted.

That a directed algebra is quasi-hereditary with 
projective standard modules follows immediately from the 
definition of quasi-hereditary algebras, see for example
\cite{DR}. This completes the proof.
\end{proof}

We would like to finish with the following easy corollary from
the above results:

\begin{corollary}\label{corollary52}
$|\mathrm{K}_n|=2|\mathrm{K}_{n-1}|+
\dim_{\mathbb{K}}(e-a_n)\mathbb{K}\mathrm{K}_n a_n$.
\end{corollary}
 
\begin{proof}
Using the proof of Proposition~\ref{proposition49} we have
\begin{multline*}
|\mathrm{K}_n|=\dim_{\mathbb{K}}\mathbb{K}\mathrm{K}_n=
\dim_{\mathbb{K}}a_n\mathbb{K}\mathrm{K}_n a_n+
\dim_{\mathbb{K}}(e-a_n)\mathbb{K}\mathrm{K}_n a_n+\\+
\dim_{\mathbb{K}}a_n\mathbb{K}\mathrm{K}_n (e-a_n)+
\dim_{\mathbb{K}}(e-a_n)\mathbb{K}\mathrm{K}_n (e-a_n)=\\
\dim_{\mathbb{K}} B +\dim_{\mathbb{K}}(e-a_n)\mathbb{K}\mathrm{K}_n a_n+
0+\dim_{\mathbb{K}} C=\\
2\dim_{\mathbb{K}}\mathbb{K}\mathrm{K}_{n-1}+
\dim_{\mathbb{K}}(e-a_n)\mathbb{K}\mathrm{K}_n a_n=
2|\mathrm{K}_{n-1}|+
\dim_{\mathbb{K}}(e-a_n)\mathbb{K}\mathrm{K}_n a_n.
\end{multline*}
\end{proof}

\vspace{0.2cm}

\noindent
G.K.: Algebra, Department of Mathematics and Mechanics, Kyiv Taras
Shevchenko University, 64 Volodymyrska st., 01033 Kyiv, UKRAINE,
e-mail: {\tt akudr\symbol{64}univ.kiev.ua}
\vspace{0.2cm}

\noindent
V.M: Department of Mathematics, Uppsala University, Box. 480,
SE-75106, Uppsala, SWEDEN, email: {\tt mazor\symbol{64}math.uu.se}

\end{document}